\documentclass[12pt]{amsart}

\DeclareMathOperator{\Br}{Br}
\DeclareMathOperator{\Gal}{Gal}

\usepackage{amscd}
\usepackage{amsfonts}
\usepackage{amssymb,latexsym,amsmath,amsfonts}

\begin{document}

\newtheorem{theorem}{Theorem}
\newtheorem{lemma}{Lemma}
\newtheorem{corollary}{Corollary}
\newtheorem{proposition}{Proposition}
\newtheorem{definition}{Definition}
\newtheorem*{example*}{Example}

\theoremstyle{definition}
\newtheorem{remark}{Remark}

\newcommand{\bs}{{\bar\sigma}}
\newcommand{\hs}{{\hat\sigma}}
\newcommand{\df}{{F^{\times}}}
\newcommand{\dk}{{K^{\times}}}
\newcommand{\F}{\mathbb{F}}
\newcommand{\Fp}{\F_{p}}
\newcommand{\Ft}{\F_{2}}
\newcommand{\Ic}{\mathcal{I}}
\newcommand{\Z}{\mathbb{Z}}
\newcommand{\lra}{\longrightarrow}
\newcommand{\gk}{\mathfrak{K}}
\newcommand{\gf}{\mathfrak{F}}

\title[Galois Module Structure of $p$th-Power Classes]
{Galois Module Structure of $p$th-Power Classes of Extensions of Degree $p$} 

\author[J\'{a}n Min\'{a}\v{c}]{J\'an Min\'a\v{c} $^{*\dagger}$}
\address{Department of Mathematics, Middlesex College, \ University of 
Western Ontario, London, Ontario \ N6A 5B7 \ CANADA}
\thanks{$^*$Research supported in part by the Natural
Sciences and Engineering Research Council of Canada and by the special Dean of
Science Fund at the University of Western Ontario.}
\thanks{$^\dag$Supported by the Mathematical Sciences Research
Institute, Berkeley}
\email{minac@uwo.ca}

\author{John Swallow}
\address{Department of Mathematics, Davidson College, Box 7046, Davidson, 
North Carolina \ 28035 \ USA}
\email{joswallow@davidson.edu}

\begin{abstract}
    For fields $F$ of characteristic not $p$ containing a primitive
    $p$th root of unity, we determine the Galois module structure of
    the group of $p$th-power classes of $K$ for all cyclic extensions $K/F$
    of degree $p$.
\end{abstract}

\maketitle

\parskip=10pt plus 2pt minus 2pt

The foundation of the study of the maximal $p$-extensions of fields $K$ containing a primitive $p$th root of unity
is a group of the $p$th-power classes of the field: by Kummer theory this group describes elementary $p$-abelian
quotients of a maximal $p$-extension. The size of the group controls the number of generators of the Galois group of
the maximal $p$-extension. Furthermore, as $K$ varies among the Galois extensions of the base field $F$, the Galois
module structure of $p$th-power classes plays a fundamental role in the investigation of the Galois group of the
maximal $p$-extension of $F$. In particular, in the nineteen-sixties Borevi\v{c} and Faddeev closely studied 
the Galois module 
structure of the $p$th-power classes of $K$ for cyclic extensions $K/F$ of degree $p$ in the case when $F$ is a local
field, and in the seventies Miki used this structure in the service of Galois embedding problems 
(see \cite{bo}, \cite{f}, and \cite{mi}). 

From a cohomological point of view, the Galois module structure of $p$th-power classes is especially important. These
Galois modules may be naturally identified with the first cohomology groups of maximal $p$-extensions with 
$\Fp$-coefficients. (See \cite[Chapter 1]{s}.) It is moreover conjectured that the entire Galois cohomology ring of 
$K$ with coefficients in $\Fp$ of degree at least one is generated by its first cohomology group. (See \cite{v} for
a proof in the case of $p=2$ as well as comments on this conjecture which was initially considered by Beilinson, 
Bloch-Kato, Lichtenbaum, and Milnor.) As a result, for the investigation of the action of the Galois group 
$\Gal(K/F)$ on the Galois cohomology of $K$ with coefficients in $\Fp$, it should be sufficient to understand the
Galois module structure of the $p$th-power classes of $K$. Furthermore this is necessary for the study of a
Galois cohomology ring of $F$ with coefficients in $\Fp$ using the Lyndon-Hochschild-Serre spectral sequence 
attached to a group extension of an absolute Galois group of $K$ by $\Gal(K/F)$. 

Hence the question of determining the Galois module structure of the
$p$th-power classes of $K$ is a fundamental question related to both the structure 
of the Galois cohomology of 
$F$ and the structure of the maximal $p$-extensions of $F$. 

It is surprising that one can determine this Galois module structure for all cyclic field
extensions of degree $p$, when the base field contains a primitive
$p$th root of unity, and that the proofs are elementary.  Just
as in the case of local fields, the structure depends upon only the norm subgroup, its
intersection with the subgroup of $p$th roots of unity, and the multiplicative subgroup of the base field. 
When we first approached these results, the 
main tools employed were Waterhouse's elegant
treatment of Galois closures of certain Kummer extensions, as well as 
a criterion for the realization of a nonabelian Galois 
group of order $p^3$ and exponent $p$ as a Galois group over our base field. We were able, however, to simplify our 
proofs considerably. In particular we require in this paper only Hilbert's Theorem 90, 
basic Kummer theory, and linear algebra over $\Fp[\Gal(K/F)]$. 
Although Galois embedding problems do not appear in 
our exposition,
our results are closely related to those in \cite{ma}, \cite{math}, and \cite{wat}, and this relationship remains 
worthy of precise formulation. This will be a topic in subsequent work.  

It is clear that one can generalize the main theorems to a larger family of field extensions $K/F$  
than we are considering here. Also, these results may be applied to  
calculations in Galois cohomology, and the study 
of Galois groups of maximal $p$-extensions of fields. We concentrate here, however, on the short exposition of a 
solution of a main significant problem. 

\section{Main Theorems}

Let $F$ be a field of characteristic not $p$ containing a primitive
$p$th root of unity $\xi_{p}$.
Let $K$ be a cyclic extension $K=F(\root{p}\of{a}),\dk=K\setminus\{0\} \mbox{ and } J=\dk/\dk^{p}$. 

Denote $\Gal(K/F)=\langle \sigma\rangle$ by $G$, where $\sigma\in\Gal(K/F)$ such that 
$\root{p}\of{a}^{\sigma-1}=\xi_{p}$.  
For $G$-modules $V$, denote by $V^{G}$ the fixed submodule of $V$. 

Maps related to the norm will play an important role, and as a result we use $N$ in several different ways. 
As an element of $\Fp[G]$ we let $N:=1+\sigma+\dots + \sigma^{p-1}$.
We also view $N$ as the resulting $G$-homomorphism $V \lra V$. Then we
write $N(J)$ for the image of $J$, and $N(\dk)$ for the image of $\dk$. 

The following two theorems classify all modules $J$ using only information about the extension $K/F$ encoded in the 
multiplicative subgroup of $F$. 

\begin{theorem}\label{th:pnot2}
    Suppose $p>2$.  As an $\Fp[G]$-module, $J$ decomposes as
    \begin{equation*}
	J=X\oplus Y\oplus Z
    \end{equation*}
    where
    \begin{enumerate}
	\item $X$ is an indecomposable $\Fp[G]$-submodule, of
	dimension $1$ if $\xi_{p}\in N(\dk)$ and of dimension $2$ if
	$\xi_{p}\not\in N(\dk)$; 
	\item $Y$ is a free $\Fp[G]$-submodule with $Y^{G}=N(J)$; and
	\item $Z$ is a trivial $\Fp[G]$-module.
    \end{enumerate}
    Moreover, the rank of any maximal free $\Fp[G]$-submodule of $J$
    is equal to the rank of $N(J)$ as an $\Fp$-module.
\end{theorem}

\begin{theorem}\label{th:p2}
    Suppose $p=2$.  As an $\Ft[G]$-module, $J$ decomposes as
    \begin{equation*}
	J=X\oplus Y\oplus Z
    \end{equation*}
    where
    \begin{enumerate}
	\item $X$ is of dimension $1$ if $-1\in N(\dk)$ and is $\{0\}$
	otherwise;
	\item $Y$ is a free $\F_{2}[G]$-submodule with $Y^{G}=N(J)$; and
	\item $Z$ is a trivial $\F_{2}[G]$-module.
    \end{enumerate}
    Moreover, the rank of any maximal free $\Ft[G]$-submodule of $J$
    is equal to the rank of $N(J)$ as an $\F_{2}$-module.
\end{theorem}

The results in Theorem~1 and Theorem~2 imply that each $\Fp[G]$-module $J$ can be written as a sum of cyclic 
$\Fp[G]$-modules. Because each cyclic $\Fp[G]$-module is an indecomposable module, and each ring of the 
$\Fp[G]$-endomorphisms of such a module is a local ring, we can conclude that our decompositions are unique. 
This follows from the Krull-Schmidt-Remak-Azumaya theorem. (See \cite[page 24]{fa}.) 

In the next theorem, which is a supplement to both Theorems~1 and 2, we determine the multiplicity of each cyclic 
$\Fp[G]$-module in a decomposition of an $\Fp[G]$-module $J$ using the arithmetical invariants associated with a field
extension $K/F$. We provide a direct proof of this theorem, and in particular we do not use the Krull-Schmidt-Remak-Azumaya
theorem quoted above. We use the same notation as above, and we also set 
$\Upsilon(K)=1$ or $0$ depending upon whether $\xi_{p}\in N(\dk)$ or $\xi_{p}\notin N(\dk)$. For each 
set $A$ we denote by $\vert A \vert$ the cardinal number of $A$. Observe that each nonzero cyclic 
$\Fp[G]$-module $M$ is characterized by its $\dim_{\Fp}M\in\{1,2,\dots,p\}$.

\begin{theorem} Let $J\cong\bigoplus_{i=1}^{p}\left(\bigoplus_{j\in \gk_{i}} M_{i,j}\right)$ be a 
decomposition of $J$ into $\Fp[G]$-cyclic modules $M_{i,j}$, where $\dim_{\Fp} M_{i,j}=i$ and 
$\gk_{i}$ are index sets for each $i\in \{1,2,\dots,p\}$.

Then we have

\begin{enumerate}
	\item $\vert\gk_{1}\vert+1=2\Upsilon(K)+\dim_{\Fp}(\df/N(\dk))$;  
	\item If $p>2$ then $\vert\gk_{2}\vert=1-\Upsilon(K)$;  
	\item[(\,i\,)] If $3\leq i\leq (p-1) \mbox{ then } \gk_{i}=\emptyset$; and 
	\item[(p)] $\vert\gk_{p}\vert=\dim_{\Fp} N(J)$.
\end{enumerate} 
\end{theorem}

The proofs of Theorems~\ref{th:pnot2}, \ref{th:p2} and~3 are presented in Section~3. 

As a consequence of these theorems and their proofs, we have the following
results on extremal cases. A short proof of each is given in the last section of this paper. 

Recall that if $p>2$ then $a\in F\setminus F^{p}$ is \emph{$p$-rigid}
if, whenever the cyclic algebra $(\frac{a,b}{F,\xi_{p}})$ is of trivial class in
$\Br(F)$, then $[b]\in \langle[a]\rangle$ in $\df/\df^p$.  For a
discussion of $p$-rigidity and $G$-invariant modules $J$, see
\cite{war}. For the basic properties of cyclic algebras and Brauer groups, see \cite[Chapters~14 and 15]{p}. 

\begin{corollary} \
    \begin{enumerate}
    \item $J$ is a free $\Fp[G]$-module precisely when $p=2$,
    $-1\not\in N(\dk)$, and $\df=N(\dk) \cup -N(\dk)$. 

    \item If $p>2$ then $J$ contains no free direct summand precisely when $a$
    is $p$-rigid. If $p=2$ then $J$ contains no free direct summand precisely when 
	$N(\dk)/F^{\times 2}\subset\langle[a]\rangle\subset\df/F^{\times 2}$. 
	(In this case it follows that $\sqrt{-1}\in\dk$.) 

    \item If $p>2$ then $J$ is $G$-invariant precisely when $a$ is $p$-rigid 
    and $\xi_{p^2}\in \dk$. If $p=2$ then $J$ is $G$-invariant precisely when $J$ contains no free summand.
    \end{enumerate}
\end{corollary}

It is worth pointing out that in the corollary above, we can replace ``free direct summand of $J$'' by 
``free submodule of $J$'', as it is well known that each free submodule of an $\Fp[G]$-module is in fact its summand. 
(See \cite[Theorem~11.2]{c}.)

Observe that if $J$ is an $\Fp[G]$-module as above, and $\varphi\colon H\lra G$ is any isomorphism, we can consider 
$J$ to be an $\Fp[H]$-module and the isomorphism type of this module is independent of the choice of $\varphi$. 
In particular, if we have two
$\Fp[G_i]$-modules $J_i,i=1,2$ with $G\cong G_i,i=1,2$, we may consider both modules $J_i$ as modules over the
same group ring $\Fp[G]$. We shall apply this remark to the following situation. 

Suppose $p>2$ and $K_{1}=F_{1}(\root{p}\of{a})$ and $K_{2}=F_{2}(\root{p}\of{a})$
are two cyclic extensions of degree $p$. We also set 
$T_i=K_{i}^{\times}/K_{i}^{\times p}, G_{i}=\Gal(K_{i}/F_{i})$ for $i=1,2$, and  
using any isomorphisms $\varphi_i\colon G\lra G_{i}$ as above
we consider both modules $T_{1}$ and $T_{2}$ as modules over $\Fp[G]$. (We write here $T_{i}$ instead of $J_{i}$
so that we avoid possible confusion with our notation for the socle series employed later in this paper.) 

Using this notation we shall formalize our classification of $\Fp[G]$-modules $J$ 
by means of arithmetic invariants associated with the multiplicative group $F^{\times}$ as follows. 

\begin{corollary} The $\Fp[G]$-modules $T_{1}$ and $T_{2}$ are isomorphic if and only if the following three conditions
are valid:

\begin{enumerate} 
\item $\xi_p\in N(K_{1}^{\times})\Leftrightarrow \xi_p\in N(K_{2}^{\times})$. 
\item $\dim_{\Fp}N(K_{1}^{\times})/F^{\times p}_{1}=\dim_{\Fp} N(K_{2}^{\times})/F^{\times p}_{2}$. 
\item $\dim_{\Fp} F^{\times}_{1}/N(K_{1}^{\times})=\dim_{\Fp}F^{\times}_{2}/N(K_{2}^{\times})$. 
\end{enumerate}

\end{corollary}

For $p=2$ we obtain a similar Corollary~3. We use the  
same notation as in the case of $p>2$ above. Moreover we set 
$\Upsilon(K_{i})=1$ or $0$ depending upon whether $-1\in N(K_{i}^{\times})$ or $-1\notin N(K_{i}^{\times})$, for both
$i=1,2$. 

\begin{corollary} The $\Ft [G]$-modules $T_{1}$ and $T_{2}$ are isomorphic if and only if the following two
conditions are valid:

\begin{enumerate}
\item $2\Upsilon(K_{1})+\dim_{\Ft}F_{1}^{\times}/N(K_{1}^{\times})=2\Upsilon(K_{2})+\dim_{\Ft}F_{2}^{\times}/N(K_{2}^{\times})$.
\item $\Upsilon(K_{2})+\dim_{\Ft}N(K_{1}^{\times})/F_{1}^{\times 2}=\Upsilon(K_{1})+\dim_{\Ft}N(K_{2}^{\times})/F_{2}^{\times 2}$.

\end{enumerate}

\end{corollary}

\section{Cyclic Modules and Fixed Elements of $J$}

This section contains a few basic facts used in our proofs of Theorems~\ref{th:pnot2} and \ref{th:p2}. 
For any subset $\{x,y,\dots\}$ of $V$ we denote by $\langle x,y,\dots \rangle$ 
the $\Fp$-submodule of $V$ spanned by $\{x,y,\dots\}$.

We freely use basic Kummer theory (see \cite[Chapter~6, Section~2]{at}), and, depending upon the context, 
we use additive and multiplicative notation for our $\Fp[G]$-modules occurring as submodules
of $\dk/\dk^{p}$. 

Let $A = \oplus_{j=0}^{p-1} \Fp \tau^{j}$ be a free $\Fp[G]$-module on
one generator $1$, where $\sigma$ acts by multiplication by $\tau$. 
For $i=1, \dots, p-1$, let $A_{i}$ denote the cyclic $\Fp[G]$-submodule
\begin{equation*}
    A_{i} = \langle(\tau-1)^{i}, (\tau-1)^{i+1}, \dots,
    (\tau-1)^{p-1}\rangle, \mbox{ and set also } A_{p}=\{0\}.
\end{equation*}
The dimension of $A/A_{i}$ is $i$, and the $\Fp[G]$-modules $A/A_{i}$
exhaust the isomorphism classes of indecomposable $\Fp[G]$-modules.
When describing an indecomposable $\Fp[G]$-module, we will use length synonymously with dimension.  
The only proper $\Fp[G]$-submodules of $A/A_{i}$ are cyclic and they are
the images of $A/A_{i}$ under $(\sigma-1)^{j}$ for $j=1, \dots, i-1$.

The socle series of $J$ will be important for the determination of indecomposable submodules of $J$.  
The socle series of $J$ is defined by $J_{1}=J^{G}$ and
$(J_{i}/J_{i-1}) = (J/J_{i-1})^{G}$ for $i>1$.  We have that
$J_{i}=\ker (\sigma-1)^{i}$, when $(\sigma-1)^{i}$ is considered as an endomorphism of $J$. 

Indecomposable submodules of $J$ may therefore be obtained as follows. We write the elements of $J$ as
$[\gamma],\gamma\in\dk$. 
For $[\gamma]\in J_{i}\setminus J_{i-1}$ let
$\gamma_{j}=\gamma^{(\sigma-1)^{j}}$ for $j=0, 1, \dots, i-1$.  Then
the $\Fp$-submodule $M_{\gamma}:=\langle[\gamma],[\gamma_{1}], \dots,
[\gamma_{i-1}]\rangle$ of $J$ is a cyclic $\Fp[G]$-submodule of
dimension $i$ isomorphic to $A/A_{i}$. 

In the proof of Theorem~\ref{th:pnot2} in the next section, we shall construct submodules $X,Y$ and $Z$ of $J$ such that $X$ 
is a cyclic submodule of $J$ of length $1$ or $2$, $Z$ is a trivial submodule of $J$, and $Y$ is a free 
$\Fp[G]$-submodule of $J$. An important part of the proof will be to show that $X,Y$ and $Z$ generate $J$ as an
$\Fp[G]$-module. In order to do this we need to construct a module $Y$ which is sufficiently large. 
The following lemma will
be used to show that we have enough free cyclic submodules of $J$ to construct our sufficiently large free
$\Fp[G]$-submodule $Y$ of $J$. 

\begin{lemma} Suppose $p>2$ and let $l$ denote the length of $M_{\gamma}$. 
	\begin{enumerate}
		\item[a)]{If $3\leq l\leq p$,  
		then there exists an element $[\alpha]\in J$ such that $\langle N([\alpha])\rangle=M_{\gamma}^{G}$.}
		\item[b)]{If $l=2$ and $\gamma$ cannot be written as $\gamma=a^{r/p}\gamma_{1}$ for some 
		$r\in\mathbb{Z}$ and $[\gamma_{1}]\in J_{1}$, then there exists an integer 
		$t\in\mathbb{Z}$ and $[\alpha]\in J$ such that 
		$\langle N([\alpha])\rangle=(M_{a^{t/p}\gamma})^{G}$.}
	\end{enumerate}
\end{lemma}

\begin{proof} Let $3\leq l\leq i\leq p$. We show by induction on $i$ that there exists an element
$\alpha_{i}\in\dk$ such that $\langle[\alpha_{i}]^{(\sigma-1)^{i-1}}\rangle=M_{\gamma}^{G}$. Then we may set
$\alpha:=\alpha_{p}$ and the proof of the first part of our lemma will be complete. If $i=l$ we set 
$\alpha_{l}=\gamma$. Assume now that $l\leq i< p$ and that our statement is true for $i$. 

Set $c=N(\alpha_{i})$. Since $[\alpha_{i}]^{(\sigma-1)^{p-1}}=[c]$ and $i<p$ we see that $[c]=[1]\in J$. Because
$c\in F^{\times}\cap K^{\times p}$ from Kummer theory, we conclude that $c=a^{s} f^{p}$ for some 
$f\in F^{\times}$ and $s\in\mathbb{Z}$. Then $N(\alpha_{i}/fa^{s/p})=1$. By Hilbert's Theorem~90 there exists 
an element $\omega\in\dk$ such that $\omega^{\sigma-1}=\alpha_{i}/fa^{s/p}$. Then 
$\omega^{(\sigma-1)^{2}}=\alpha_{i}^{(\sigma-1)}/\xi_{p}^{s}$. Therefore 
$\langle\omega^{(\sigma-1)^{i}}\rangle=M_{\gamma}^{G}$ and we can set $\alpha_{i+1}:=\omega$. 

Now we shall prove the second part of our lemma. Assume $l=2=i$. Proceeding in the same way as above, we see that
for $\alpha_{2}=\gamma$ we have $N(\alpha_{2})=c=a^{s}f^{p}$ for some $f\in F^{\times}$ and $s\in\mathbb{Z}$. 

As before we see that there exists an element $\omega\in\dk$ such that 
$\omega^{\sigma-1}=\alpha_{2}/fa^{s/p}$. Then $\omega^{(\sigma-1)^{2}}=(\alpha_{2}a^{-s/p})^{(\sigma-1)}=
(\gamma a^{-s/p})^{(\sigma-1)}$. Observe that $\gamma a^{-s/p}\notin J_{1}$ by our hypothesis, and therefore
the length of $M_{\gamma a^{-s/p}}$ is again $2$. 
Hence we can set $\alpha_{3}:=\omega$. We can then continue by induction on $i$ as above, concluding 
that there exists an element $\alpha:=\alpha_{p}$ such that $\langle N([\alpha])\rangle=
(M_{\gamma a^{-s/p}})^{G}$ as required.
\end{proof}

\begin{remark} In the case of $\xi_{p}\in N(\dk)$ we do not need the adjusting factor 
$a^{t/p}$ in part b) of Lemma~1. Also in this case no restriction on $\gamma$ is necessary and $\gamma$ can be any 
element of $\dk$ such that the length of $M_{\gamma}$ is $2$. 

Indeed let $\beta\in\dk$ such that $N(\beta)=\xi_{p}$. As before we have 
$\omega^{(\sigma-1)}=\alpha_{2}/fa^{s/p}$ for some $\omega\in\dk$ and $s\in\Z$. Then set 
$\theta=(\beta^{(\sigma-1)^{p-3}})^{s} \;\omega$. Then $\theta^{(\sigma-1)^{2}}=\alpha_{2}^{(\sigma-1)}$. Hence we may 
set $\alpha_{3}:=\theta$. We may then continue by induction on $i$ to conclude that there exists an element
$\alpha:=\alpha_{p}$ such that $\langle N([\alpha])\rangle =M_{\gamma}^{G}$. 
\end{remark}

In the next lemma we determine $J_{1}$ in terms of arithmetic invariants of $K/F$ encoded in the multiplicative 
group of $F$. We write $\epsilon(\df)$ for the subgroup $\df\dk^{p}/\dk^{p}$ of $J$.

\begin{lemma} \
	\begin{enumerate}
		\item[a)]{If $\xi_{p}\notin N(\dk)$ then $J_{1}=\epsilon(F^{\times})$.}
		\item[b)]{If $\xi_{p}\in N(\dk)$ then $J_{1}\cong\epsilon(F^{\times})
		\oplus\langle\delta\rangle$ where $\delta\in\dk,
		\sigma(\delta)/\delta=\lambda^{p}$ and $N(\lambda)=\xi_{p}$.}
	\end{enumerate}
\end{lemma} 

\begin{proof} Suppose that $\theta\in\dk$ such that $[\theta]\in J_{1}$. Then $\sigma(\theta)/\theta=\lambda^{p}$
for some $\lambda\in\dk$, and hence 
$N(\lambda)^{p}=1$. Hence we see that $N(\lambda)=\xi_{p}^{c}$ for some
$c\in\Z$. Now consider the case $\xi_{p}\notin N(\dk)$. Then $N(\lambda)=1$, because otherwise $\xi_{p}$ would be 
a norm of a suitable power of $\lambda$. Therefore from Hilbert's Theorem~90 we see that 
$\sigma(\theta)/\theta=\sigma(k^{p})/k^{p}$ for some $k\in\dk$. We conclude that 
$\theta/k^{p}\in\df$ and hence 
$[\theta]=[f]$ for some $f\in F^{\times}$. Therefore if $\xi_{p}\notin N(\dk)$ then 
$J_{1}=\epsilon(F^{\times})$ as required. 

Now assume that $\xi_{p}\in N(\dk)$. Then $\xi_{p}=N(\lambda)$ for some $\lambda\in\dk$ and by Hilbert's Theorem~90
there exists an element $\delta\in\dk$ such that $\sigma(\delta)/\delta=\lambda^{p}$. Then the $\Fp[G]$-submodule of
$J$ generated by $\delta$ and $\epsilon(F^{\times})$ is isomorphic with $\epsilon(F^{\times})\oplus\langle\delta\rangle$.
Now for each $[\theta]\in J_{1}, \sigma(\theta)/\theta=\nu^{p},N(\nu)=\xi_{p}^{c},c\in\Z$, and therefore we have  
$[\theta]\in\epsilon(F^{\times})[\delta]^{c}$. Hence $J_{1}\cong\epsilon(F^{\times})\oplus\langle\delta\rangle$ as
required. 
\end{proof}

\begin{remark} It is worthwhile to observe that we have the following short exact sequence:
$$0\lra\langle[a]\rangle\overset{i}{\lra}F^{\times}/F^{\times p}\overset{\epsilon}{\lra} J_{1}\overset{N}{\lra}\langle[a]\rangle,$$
\noindent where $\langle[a]\rangle$ is the subgroup of $F^{\times}/F^{\times p}$ generated by 
$[a]\in F^{\times}/F^{\times p}$,
$i$ is the inclusion map, $\epsilon$ is the natural homomorphism induced by the inclusion map 
$\df\lra\dk$ and $N$ is the map induced by the norm map $N\colon\dk\lra\df$. Moreover the map 
$N\colon J_{1}\lra\langle[a]\rangle$ is surjective if and only if $\xi_{p}\in N(\dk)$. 

The fact that $[N(\theta)]\in\langle[a]\rangle$ for each $[\theta]\in J_{1}$ follows from the fact that 
$N(\theta)\in\df\cap\dk^{p}$. The exactness at $\df/F^{\times p}$ 
follows from \linebreak Kummer theory. The exactness at $J_{1}$ can be seen as follows. 

Suppose that $[\theta]\in J_{1}$. Then $\root{p}\of{N(\theta)}\in\dk$ and $[N(\theta)]=[1]\in\df/\df^{p}$ if and only if 
$\sigma(\root{p}\of{N(\theta)})/\root{p}\of{N(\theta)}=1$. Set $\sigma(\theta)/\theta=\lambda^{p},
\lambda\in\dk$. Then the condition $\sigma(\root{p}\of{N(\theta)})/\root{p}\of{N(\theta)}=1$ translates as 
$N(\lambda)=1$. Using Hilbert's Theorem~90 as in the proof of Lemma~2, part a) above, we see that 
$N(\lambda)=1$ if and only if $[\theta]\in\epsilon(\df)$. 
  
Finally if $\xi_{p}\notin N(\dk)$ we know that $\epsilon\colon\df/F^{\times p}\lra J_{1}$ 
is surjective and hence $N\colon J_{1}\lra\langle[a]\rangle$ is trivial. If $\xi_{p}\in N(\dk)$, then 
$\epsilon\colon\df/F^{\times p}\lra J_{1}$ is not surjective and therefore $N$ must be surjective. 
\end{remark}

\section{Proofs}

\begin{proof}[Proof of Theorem~\ref{th:pnot2}]
	In the first part of this proof we shall construct some submodules $X,Y$ and $Z$ of $J$, and we shall
	show that they generate a submodule of $J$ isomorphic to $X\oplus Y\oplus Z$. In the second part of this proof
	we shall show that $X,Y$ and $Z$ generate the full module $J$. 

    	We first construct $X$.  If $\xi_{p}\in N(\dk)$, then by Remark~2 at the end of Section~2 we see that there exists
	an element $\delta\in\dk$ such that $N(\delta)=a$ and $[\delta]\in J_{1}$. Then we set 
	$X=\langle[\delta]\rangle$.  
        
    If $\xi_{p}\not\in N(\dk)$, then let $\delta=\root{p}\of{a}$ and
    $X=\langle[\delta], [\xi_{p}]\rangle \subset J_{2}$.  If
    $[\xi_{p}]=[1]$ then a root of unity $\xi_{p^{2}}$ of order
    $p^{2}$ lies in $\dk$ and we can pick $\xi_{p^{2}}$ such that $\xi_{p}=N(\xi_{p^{2}})$, contrary
    to hypothesis.  Because $[\root{p}\of{a}]^{\sigma-1}=[\xi_{p}]$,
    $X$ is isomorphic to $A/A_{2}$.  
    
    We proceed to construct $Y$.  Let $\Ic$ be an $\Fp$-basis for
    $N(J)$.  For each $[x]\in \Ic$ we construct a free $\Fp[G]$-module
    $M(x)$, as follows.  Choose a representative $x\in \df$ for $[x]$, such that  
    $x\in N(\dk)$.  Choose $\gamma\in \dk$ such that $x=N(\gamma)$. Finally set $M(x)=M_{\gamma}$.  
    
    We claim that the $M(x)$, $[x]\in \Ic$, are independent.  First we
    show by induction on the number of modules that a finite set of
    modules $M(x)$ is independent.  The base case is trivial.  Now
    suppose that $W=M(x')\cap \sum_{[x']\neq [x]} M(x)\neq \{0\}$. 
    Then $W$ contains the $1$-dimensional submodule $V=M(x')\cap J_{1}$.  If
    $M(x')=M_{\omega}$, then $V=\langle[N(\omega)]\rangle = 
    \langle[x']\rangle$.  Since $[x']\in W$, $[x']$ is a finite sum
    $\sum_{[x']\neq [x]} [m(x)]$ where $[m(x)]\in M(x)$.  By induction
    the modules $M(x)$ appearing in the sum are independent.  Then
    since $[x']$ is in $J_{1}$, each $[m(x)]\in J_{1}$ as well.  But
    then each $[m(x)]\in \langle[x]\rangle$, hence $[m(x)]\in N(J)$. 
    Since $[x']$ and the finite number of elements $[x]$ considered above are distinct elements of an
    $\Fp$-base for $N(J)$, they are independent, and we have a
    contradiction.
    
    Now if $W=M(x')\cap \sum_{[x']\neq [x]} M(x)\neq \{0\}$ where the
    sum is infinite, the same argument holds, since $[x']\in V=W\cap
    J_{1}$ is a finite sum of elements $[m(x)]$.  Hence the $M(x)$,
    $[x]\in \Ic$, are independent.
    
    Let $Y=\oplus_{\Ic} M(x)$.  Then $Y$ is a free $\Fp[G]$-module
    with a generating set in one-to-one correspondence with a
    generating set for $N(J)$ as an $\Fp$-module.  Moreover, $Y^{G} =
    \oplus M(x)^{G} = \oplus \langle[x]\rangle = N(J).$
    
    Let $Z$ be any complement in $\epsilon(\df)$ of the
    $\Fp$-submodule of $J$ generated by $N(J)$ and $X\cap
    \epsilon(\df)$.  Clearly $Z$ is a trivial $\Fp[G]$-module. 

    We claim that $X$, $Y$, and $Z$ are independent $\Fp[G]$-submodules of $J$. If $W=Y\cap
    Z\neq \{0\}$, then $W$ is a submodule of $\epsilon(\df)$.  Let
    $[x]\in W$.  Then $[x]$ lies in $Y\cap J_{1}=N(J)$, as well as in
    $Z$.  But $Z$ is a complement of $\epsilon(\df)$ of a submodule
    containing $N(J)$, a contradiction. Therefore the $\Fp[G]$-submodule of $J$ generated by $Y$ and $Z$ is isomorphic with
    $Y\oplus Z$. 	
    
    Now suppose $W=X\cap (Y\oplus Z)\neq \{0\}$. Considering $W\cap J_{1}$ as above, we find that $X$ is of
    dimension $2$, $\xi_{p}\not\in N(\dk)$, and $W=\langle[\xi_{p}]\rangle$.  If $[\xi_{p}]\in
    N(J)$ then $\xi_{p}a^{c}\in N(\dk)$ for some $c\in \mathbb{Z}$.  However,
    since $N(\root{p}\of{a})=a$, we have $\xi_{p}\in N(\dk)$, a
    contradiction.  Therefore $W\cap Y=\{0\}$ and $[\xi_{p}]= y + z$,
    $y\in Y$, $z\in Z$, with $z\neq 0$.     
    Hence $y\in J_{1}\cap Y=N(J)$, and 
    $z=[\xi_{p}]-y$ is the relation in $\epsilon(\df)$.  But $Z$ is a
    complement in $\epsilon(\df)$ of the submodule generated by $N(J)$
    and $(X\cap \epsilon(\df))=\langle[\xi_{p}]\rangle$, a
    contradiction.
    
    Let $\tilde J$ be the $\Fp[G]$-submodule of $J$ generated by $X\cup Y\cup Z$. As we have shown above, 
	$\tilde J\cong X\oplus Y\oplus Z$.  
    We show that $J=\tilde J$ by induction on the socle series of $J$.  
    
    First we show $J_{1}\subset \tilde J$. From our construction of modules $X,Y$ and $Z$ we see that 
	$\epsilon(\df)\subset\tilde{J}$. Therefore from Lemma~2a) in Section~2 we see that if $\xi_{p}\notin N(\dk)$ then
	$J_{1}\subset\tilde{J}$. If $\xi_{p}\in N(\dk)$ then from the definition of $X$ and Lemma~2b), we see that
	again $J_{1}\subset\tilde{J}$. 
    
    For the inductive step, assume $J_{i}\subset \tilde J, 1\leq i< p$, and let 
    $[\gamma]\in J_{i+1}\setminus J_{i}$. 

	Let $\langle[b]\rangle=M_{\gamma}^{G}$. From our Lemma~1a) in Section~2 we see that if $i\geq 2$ then 
    	$[b]\in N(J)$. If $i=1$ and $\xi_{p}\in N(\dk)$ we see from Remark~1 in Section~2 that again $[b]\in N(J)$. 
	If $\xi_{p}\notin N(\dk)$ then $[\root{p}\of{a}]\in X$ and therefore it is sufficient to show that
	$[a^{t/p}\gamma]\in\tilde{J}$ for some $t\in\Z$. Therefore in this case if $\gamma=a^{r/p}\gamma_{1}$, where
	$r\in\mathbb{Z}$ and $[\gamma_{1}]\in J_{1}$, we see that $[\gamma]\in\tilde{J}$ and if $\gamma$ cannot be written
	in this form, we replace $\gamma$ by a suitable 
	$a^{t/p}\gamma$, and we use Lemma~1b) to conclude that $(M_{a^{t/p}\gamma})^{G}\subset N(J)$. 
	Thus using Lemma~1 we have reduced our 
	investigation to the case when $[b]\in N(J)$. 
	
	We may write $[b]=\sum_{\mathcal{I}}c_{x}[x]$  
    with $c_x\in \Fp$ and almost all $c_x=0$.  Now for each $[x]$,
    $M(x)=M_{\omega}$ for some $\omega\in\dk$ with $N([\omega])=[x]$. 
    For each $[x]$, set $\gamma(x) = \omega^{(\sigma-1)^{p-i-1}}$, so
    that $[\gamma(x)] \in J_{i+1}$ and $[\gamma(x)]^{(\sigma-1)^i} =
    [x]$.  Now let $[\Gamma]=\sum_{\mathcal{I}}c_{x}[\gamma(x)]$. The 
    application of $(\sigma-1)^i$ to $[\Gamma]$ results in $\sum_{\mathcal{I}}c_{x}[x]=[b]$.
    Hence $[\gamma/\Gamma]$ lies in the kernel of
    $(\sigma-1)^i$, and $[\gamma/\Gamma]\in J_i\subset \tilde J$ by
    induction.  Since $[\Gamma]\in \tilde J$, $[\gamma]\in \tilde J$.
\end{proof}

The proof of Theorem~2 below follows the previous proof with some modifications. Since in the case of 
$p=2$ the length of the socle series of $J$ is at most $2$, the proof is simpler. Moreover, we make use of the
well-known fact on elements of identical norm in quadratic extensions:

If $K=F(\sqrt{a}),a\in F^{\times}\setminus F^{\times 2},\gamma_{1},\gamma_{2}\in K^{\times}$ and 
$N(\gamma_{1})=N(\gamma_{2})$ then $[\gamma_{1}]\in\epsilon(F^{\times})[\gamma_{2}]$. Indeed by Hilbert's Theorem~90
we have $[\gamma_{1}]=[N(k)][\gamma_{2}]$ for some $k\in\dk$. (See \cite[page~202]{l}.)

\begin{proof}[Proof of Theorem~2]
If $-1\in N(K^{\times})$, then we set $X=\langle[\delta]\rangle\subset J_{1}$ with $\delta$ such that 
$\delta^{\sigma-1}=\beta^{2}$ and $N(\beta)=-1$ (observe that then $[N(\delta)]=[a]\in F^{\times}/F^{\times 2})$. If 
$-1\notin N(K^{\times})$ then we set $X=\{0\}$. 

Choose $\Ic$ to be an $\F_{2}$-basis of $N(J)$. For each $[x]\in\Ic$ choose a representative 
$x\in F^{\times}$ of $[x]$ and $\omega\in K^{\times}$ such that $x=N(\omega)$. Then set 
$M(x)=M_{\omega}$ and let $Y$ be the submodule of $J$ generated by all $M(x),x\in\Ic$. The same argument as in the
proof of Theorem~1 shows that $Y$ is a free $\F_{2}[G]$-submodule of $J$, and that 
$Y\cong\oplus M(x)$. Finally let $Z$ be any complement in $\epsilon(F^{\times})$ of the $\F_{2}$-module
$N(J)$. As before we check that $X,Y$ and $Z$ are independent, $Z$ is a trivial 
$\F_{2}[G]$-module, and we let $\tilde J\cong X\oplus Y\oplus Z$ be the $\F_{2}[G]$-submodule of $J$ generated by 
$X,Y$ and $Z$. 

Assume now that $[\gamma]\in J_{1}$. Then $N(\gamma)\in F^{\times}\cap K^{\times 2}=F^{\times 2}\cup aF^{\times 2}$. 
If $N(\gamma)\in F^{\times 2}$ then $[\gamma]\in\epsilon(F^{\times})\subset\tilde{J}$. 
If $N(\gamma)\in aF^{\times 2}$ then 
$[\gamma]\in[\delta]\epsilon(F^{\times})\subset\tilde{J}$. Thus $J_{1}\subset\tilde{J}$. Now suppose that 
$[\gamma]\in J\setminus J_{1}$. Let $\lambda=\gamma^{\sigma-1}$ and since 
$\sigma-1=\sigma+1$ in $\F_2[G], \mbox{ we see that } [\lambda]=[N(\gamma)]\neq[1]$. 
Hence for a suitable $\Gamma\in Y$ we have $[\gamma/\Gamma]\subset J_1$.  
Therefore $\tilde{J}=J$ as required.
\end{proof}

\begin{proof}[Proof of Theorem~3]
	We shall first show that if $3\leq i\leq p-1$ then $\mathfrak{K}_{i}=\emptyset$. Suppose that this is not true 
	and that there exists $M_{i,k}= M,3\leq i\leq(p-1), k\in\mathfrak{K}_{i}$, appearing as a summand in our 
	decomposition of $J$. 

	Let $m$ be a generator of $M$. Using the decomposition $J=X\oplus Y\oplus Z$ 
	described in Theorem~1, we can write $m=u+v$, where $u\in Y$ and $v\in X\oplus Z$. Using the fact that
	$v\in J_{2}$ and $3\leq i$ we see that we can find an element $n\in M\setminus\{0\}$ and 
	$y\in Y\setminus Y\cap J_{p-1}$ such that $n=(\sigma-1)^{k}y$ for a suitable $k\in\mathbb{N}$. 

	Now write $y=\sum_{s=1}^{t}m_{s}$ where $m_{s}\in M_{a,j},1\leq a\leq p$ and $j\in\mathfrak{K}_{a}$. 
	Of course we assume that
	we choose at most one element from each $M_{a,j}$. Then from the equality 
	$n=\sum_{s=1}^{t}(\sigma-1)^{k}m_{s}$ we see that we may assume that $n=(\sigma-1)^{k}m_{1}$ and 
	$(\sigma-1)^{k}m_{s}=0$ for $s\in\{2,\dots,t\}$. Further from our equalities 
	$y=\sum_{s=1}^{t}m_{s},(\sigma-1)^{k}m_{s}=0$ for $s\in\{2,\dots,t\}$ and $y\in J_{p}\setminus J_{p-1}$ 
	we see that $m_{1}\in J_{p}\setminus J_{p-1}$. But this implies that $M$ is a free summand of $J$ - a
	contradiction. Hence if $3\leq i\leq p-1$ then $\mathfrak{K}_{i}=\emptyset$ as we claimed. 

	Now we consider the case of $i=p$. Since $N(J)=\oplus_{j\in\mathfrak{K}_{p}} N(M_{p,j})$ and each 
	$N(M_{p,j})$ is a $1$-dimensional $\Fp$-module, we see that $\vert\gk_{p}\vert=\dim_{\Fp} N(J)$. 
	
	Let us now consider the case of $i=2$ and $p>2$. We have 
	$$\vert\gk_{2}\vert=\dim_{\Fp}(((\sigma-1)J)^{G}/N(J)).$$ 
	Therefore we see that 
	$\vert\mathfrak{K}_{2}\vert$ is an invariant of $J$ which does not depend upon any particular choice of the
	decomposition of $J$ into a sum of cyclic submodules of $J$. Therefore in order to determine 
	$\vert\mathfrak{K}_{2}\vert$ we may use our decomposition $J=X\oplus Y\oplus Z$, where $X,Y$ and $Z$ are each
	thought of as a sum of their cyclic submodules. Thus in the case of $p>2$, from Theorem~1 we conclude that 
	$\vert\mathfrak{K}_{2}\vert=1-\Upsilon(K)$. 
	
	Finally observe that if $p>2$ then $\vert\gk_{1}\vert+\vert\gk_{2}\vert=\dim_{\Fp}J_{1}/N(J)$, and if 
	$p=2$ then $\vert\gk_{1}\vert=\dim_{\F_{2}}J_{1}/N(J)$. 
	Applying Lemma~2 and observing that if $p>2$ then $a\in N(\dk)$ and if 
	$p=2$ then $a\in N(\dk)$ if and only if $-1\in N(\dk)$, we obtain 
	$$1+\vert\gk_{1}\vert=\dim_{\Fp}(\df/N(\dk))+2\Upsilon(K).$$
\end{proof}

\begin{proof}[Proof of Corollary~1] \

     	(1): From Theorem~3 we see that if $p>2$ then $J$ is a free $\Fp[G]$-module 
	if and only if $\gk_{1}=\emptyset=\gk_{2}$, and if $p=2$ then $J$ is a free $\Fp[G]$-module if and only if
	$\vert\gk_{1}\vert=\emptyset$. 
	Therefore if $J$ is a free $\Fp[G]$-module, $p=2$ and $-1\notin N(\dk)$. 
	Moreover assuming $p=2,-1\notin N(\dk)$ we see that
	$\gk_{1}=\emptyset$ if and only if $\df=N(\dk)\cup - N(\dk)$.
     
     (2): $J$ contains no free direct summand precisely when $N(J)$ is
     trivial. If $p>2$ then $N(J)$ is trivial precisely when $a$ is $p$-rigid. Suppose now that $p=2$. 
	Then we see that $N(J)$ is trivial if and only if $N(\dk)/F^{\times 2}
	\subset\langle[a]\rangle\subset\df/F^{\times 2}$. Since 
	$[-1]=[-a]\in N(J)$ we see that in this case $\sqrt{-1}\in\dk$. 
	
     (3): Suppose first that $p>2$. If $J$ is $G$-invariant then $J$ contains no free direct
     summand.  Hence $a$ is $p$-rigid.  If $a$ is $p$-rigid then
     $N(J)$ is trivial, and $J$ is $G$-invariant if
     additionally $\xi_{p}\in N(\dk)$.  By the definition of
     $p$-rigidity, $[\xi_p]\in \langle[a]\rangle$ in $\df/\df^{p}$,
     which is equivalent with $\xi_{p^{2}}\in \dk$. Now if $p=2$ then $J$ is $G$-invariant precisely when $J$ contains
	no free direct summand.
\end{proof}

\begin{proof}[Proof of Corollary~2]
Suppose first that our three conditions mentioned in this corollary are valid. Then we see that 
if we use Theorem~\ref{th:pnot2} and
its proof, we can write $T_{i}=K^{\times}_{i}/K^{\times p}_{i}$ and $T_{i}=X_{i}\oplus Y_{i}\oplus Z_{i},i=1,2$, where
modules $X_{i},Y_{i}$ and $Z_{i}$ are described in Theorem~1 and its proof.   

From condition~1 of our corollary we see that $X_{1}\cong X_{2}$, from condition~2 we see that 
$Y_{1}\cong Y_{2}$, and finally with condition~3 together with our choice of $Z_{i}$ as an $\Fp$-complement of 
the $\Fp$-submodule of $T_{i}$ generated by $X_{i}\cap\epsilon
(F^{\times}_{i})$ and $N(T_{i})$ in $\epsilon(F^{\times}_{i})$ 
we see that $Z_{1}\cong Z_{2}$. Hence if conditions (1), (2) and (3) are valid, then $T_{1}\cong T_{2}$. 

Now assume that $T_{1}=X_{1}\oplus Y_{1}\oplus Z_{1}$ where submodules $X_{1},Y_{1},Z_{1}$ of $T_{1}$ are constructed 
as in the proof of Theorem~1 and $\psi$ is an $\Fp[G]$-isomorphism $T_{1}\lra T_{2}$. Set
$X_{2}=\psi(X_{1}),Y_{2}=\psi(Y_{1})$ and $Z_{2}=\psi(Z_{1})$. Then $T_{2}=X_{2}\oplus Y_{2}\oplus Z_{2}$. We also see 
that $N(T_{2})=N(Y_{2})\cong N(Y_{1})= N(T_{1})$. Therefore we have 
$\dim_{\Fp} N(K_{1}^{\times})/F^{\times p}_{1}=\dim_{\Fp} N(K_{2}^{\times})/F^{\times p}_{2}$.  

We shall now assume that $X_{1},Y_{1}$ and $Z_{1}$ are decomposed into their cyclic summands, and we shall apply 
Theorem~3. 

We obtain $\xi_{p}\in N(K_{1}^{\times})\Leftrightarrow \gk_{2}\neq\emptyset\Leftrightarrow\xi_{p}\in N(K_{2}^{\times})$. 
Finally from Theorem~3 we see that $\dim_{\Fp}(F_{1}^{\times}/N(K_{1}^{\times}))=
\dim_{\Fp}(F_{2}^{\times}/N(K_{2}^{\times}))$ as required. 
\end{proof}
Finally we shall prove Corollary~3. Observe that this proof is a straightforward corollary of Theorem~3. 
We write $\gk_{2,1}$ for the index set $\gk_{2}$ used in Theorem~3 when applied to the field extension $K_{1}/F_{1}$. 
Similarly we extend this notation for other index sets under consideration. 

\begin{proof}[Proof of Corollary~3]
From Theorem~3 we see that $T_{1}\cong T_{2}$ if and only if $\vert\gk_{1,1}\vert=\vert\gk_{1,2}\vert$ and
$\vert\gk_{2,1}\vert=\vert\gk_{2,2}\vert$. Rewriting these equalities we obtain $T_{1}\cong T_{2}$ if and only if 
we have the following two equalities. 

\begin{enumerate}
	\item $2\Upsilon(K_{1})+\dim_{\Ft}F_{1}^{\times}/N(K_{1}^{\times})=2\Upsilon(K_{2})+\dim_{\Ft}F_{2}^{\times}/N(K_{2}^{\times})$,
	\item $\Upsilon(K_{2})+\dim_{\Ft}N(K_{1}^{\times})/F_{1}^{\times 2}=\Upsilon(K_{1})+\dim_{\Ft}N(K_{2}^{\times})/F_{2}^{\times 2}.$
\end{enumerate}
\end{proof}
\section{Acknowledgments}

We gratefully acknowledge the support of the Mathematical Sciences Research Institute in Berkeley, California, 
where during the fall of 1999 the authors met, for their subsequent discussions led to the results in this paper. 
The first author is also very appreciative of the stimulating discussions with Adrian Wadsworth.

\end{document}